\documentclass[12pt]{amsart}
\usepackage{amscd,amssymb}

\usepackage[graph,frame,poly,arc]{xy}  
\usepackage[plainpages,backref,urlcolor=blue]{hyperref}

\theoremstyle{plain}
\newtheorem{thm}[subsection]{Theorem}

\newtheorem{cor}[subsection]{Corollary}

\theoremstyle{definition}
\newtheorem{rk}[subsection]{Remark}

\newtheorem{ex}[subsection]{Example}

\numberwithin{equation}{section}
\newcommand{\A}{{\mathcal A}}

\newcommand{\B}{{\mathcal B}}

\newcommand{\al}{{\alpha}}

\newcommand{\C}{\mathbb{C}}

\newcommand{\PP}{\mathbb{P}}

\DeclareMathOperator{\im}{im}

\begin{document}

\title [Some remarks on plane curves related to freeness ]
{Some remarks on plane curves related to freeness }

\author[Alexandru Dimca]{Alexandru Dimca$^1$}
\address{Universit\'e C\^ ote d'Azur, CNRS, LJAD, France and Simion Stoilow Institute of Mathematics,
P.O. Box 1-764, RO-014700 Bucharest, Romania}
\email{Alexandru.DIMCA@univ-cotedazur.fr}

\thanks{$^1$ partial support from the project ``Singularities and Applications'' - CF 132/31.07.2023 funded by the European Union - NextGenerationEU - through Romania's National Recovery and Resilience Plan.}

\subjclass[2010]{Primary 14H50; Secondary  13D02}

\keywords{plane curve, free curve, Milnor algebra, minimal resolution, Tjurina number, Betti numbers, line arrangement}

\begin{abstract} Let $C$ be a reduced complex projective plane curve,
and let $d_1$ and $d_2$ be the first two smallest exponents of $C$.
 For a free curve $C$ of degree $d$, there is a simple formula relating $d,d_1, d_2$ and the total Tjurina number of $C$. Our first result discusses how this result changes when the curve $C$ is no longer free. For a free line arrangement, the Poincaré polynomial coincides with the Betti polynomial $B(t)$ and with the product $P(t)=(1+d_1t)(1+d_2t)$. Our second result shows that for any curve $C$, the difference $P(t)-B(t)$ is a polynomial $a t +bt^2$, with $a$ and $b$ non-negative integers. Moreover $a =0$ or $b=0$  if and only if  $C$ is a free line arrangement. 
Finally we give new bounds for the second exponent $d_2$ of a line arrangement $\A$, the corresponding lower bound being an improvement of a result by H. Schenck concerning the relation between the maximal exponent of $\A$ and the maximal multiplicity of points in $\A$.
\end{abstract}
 
\maketitle


\section{Introduction} 

We denote by $S=\C[x,y,z]$  the polynomial ring in three variables $x,y,z$ with complex coefficients, and by $C:f=0$  a reduced curve of degree $d\geq 3$ in the complex projective plane $\PP^2$. 
Let $J_f$ be the Jacobian ideal of $f$, i.e., the homogeneous ideal in $S$ spanned by the partial derivatives $f_x,f_y,f_z$ of $f$, and  by $M(f)=S/J_f$ the corresponding graded quotient ring, called the Jacobian (or Milnor) algebra of $f$.
Consider the graded $S$-module of Jacobian syzygies of $C$ or, equivalently, the module of derivations killing $f$, namely
\begin{equation}
\label{eqD0}
D_0(C)= \{\theta \in Der(S) \ : \ \theta(f)=0\}.
\end{equation}
In the sequel, a derivation $\theta=a\partial_x+b\partial_y +c \partial_z$
will be identified to the associated triple $\rho=(a,b,c) \in S^3$.
We say that $C$ is an {\it $m$-syzygy curve} if  the module $D_0(C)$ is minimally generated by $m$ homogeneous syzygies, say $\rho_1,\rho_2,\ldots ,\rho_m$, of degrees $d_j=\deg \rho_j$ ordered such that $$d_1\leq d_2 \leq \ldots \leq d_m.$$ 
We call these degrees $(d_1, \ldots, d_m)$ the {\it exponents} of the curve $C$. 
The smallest degree $d_1$ is sometimes denoted by ${\rm mdr}(f)$ and is called the minimal degree of a Jacobian relation for $f$. 

The $S$-module $D_0(C)$ is free if and only if $m=2$, and in this case the curve $C$ is said to be free. The exponents $(d_1,d_2)$ of a free curve $C$ of degree $d$ are known to satisfy two key relations, namely
\begin{equation}
\label{eqF1}
d_1+d_2=d-1
\end{equation}
and
\begin{equation}
\label{eqF2}
d_1d_2= (d-1)^2-\tau(C),
\end{equation}
where  $\tau(C)$ denotes the total Tjurina number of $C$, see for instance \cite{Dmax,dPCTC}. It is known that the condition
\eqref{eqF1}
characterizes the free curves, that is a curve $C$ which is not free satisfies
\begin{equation}
\label{eqF3}
d_1+d_2 \geq d,
\end{equation}
with equality exactly for the plus-one generated curves, see \cite[Theorem 2.3]{minTjurina} for details. 

For an arbitrary curve $C$, there are many results involving the first exponent $d_1$, see for instance \cite{dPCTC}, with additional information when $C$ is a line arrangement, see for instance \cite{Mich}.
{\it The main theme of this paper is to give some information involving the first two exponents $d_1$ and $d_2$ of an arbitrary curve $C$.}
Except from \eqref{eqF3} we are not aware of any such results in the existing literature.

The first main result of this note is to completely characterize the curves satisfying the second condition \eqref{eqF2}, see Theorem \ref{thm2}. Beyond the free curves, the equality \eqref{eqF2} is satisfied only for 3-syzygy curves with $d_3=d-1$ and for all the other curves one has
\begin{equation}
\label{eqF21}
d_1d_2 > (d-1)^2-\tau(C).
\end{equation}

When $C=\A$ is a line arrangement in $\PP^2$,
one can define its Poincaré polynomial $\pi(\A,t)$ and  it is well known that
$$\pi(\A,t)=(1+d_1t)(1+d_2t)$$
if $\A$ is free with exponents $(d_1,d_2)$, see \cite{OT,DHA}. Moreover, the Poincaré polynomial $\pi(\A,t)$ coincides with the Betti polynomial $B(M(\A))(t)$ of the complement $M(\A)=\PP^2 \setminus \A$, see \eqref{eqE1} for the general definition of the Betti polynomial and  see \cite{OT,DHA} for the equality $\pi(\A,t)=B(M(\A))(t)$.
Hence the above relation may be written as
\begin{equation}
\label{eqF4}
B(M(\A))(t)=(1+d_1t)(1+d_2t).
\end{equation}
The second main result of this note is Theorem \ref{thm3}, which describes the relation between $B(M(C))(t)$ and the product $(1+d_1t)(1+d_2t)$ for any  curve $C$ in $\PP^2$. It shows that the difference
$$(1+d_1t)(1+d_2t)-B(M(C))(t)$$
is a polynomial $a(C) t +b(C) t^2$, with $a(C)$ and $b(C)$ non-negative integers. Moreover $a(C) =0$ or $b(C)=0$ implies that $C$ is a free line arrangement.

As applications of Theorems \ref{thm2} and \ref{thm3} we show that any irreducible component of a free curve with $d_1=1$ is rational, see Remark \ref{rk1} and Corollary \ref{cor31}. 

 The first result in Section 5 gives a combinatorial restriction on the number of multiple points on any line in the arrangement $\A$, which a potential minimal counter-example $\A$ to Terao's Conjecture must satisfy, see Theorem \ref{thm4}. This result was already known, see Remark \ref{rk40}, but we hope our approach may be useful for some readers.
 The main result in Section 5  gives new bounds for the second exponent $d_2$ of any line arrangement $\A$,  see Theorem \ref{thm5}.
 The corresponding lower bound is an improvement of a claim by H. Schenck, see  \cite[Lemma 5.2]{Sch}, concerning the relation between the maximal exponent of $\A$ and the maximal multiplicity of points in $\A$.
 
 \medskip
 
 We thank Lukas K\"uhne  and Piotr Pokora for very useful discussions related to
 Theorems \ref{thm4} and \ref{thm5}.

\section{On plane curves satisfying $\tau(C)=(d-1)^2-d_1d_2$} 

We recall first the construction of the Bourbaki ideal $B(C,\rho_1)$ associated to a degree $d$ reduced curve $C:f=0$ and to a minimal degree non-zero syzygy $\rho_1 \in D_0(C)$, see \cite{DStJump}
as well as \cite{JNS} for a recent, more complete approach.
For any choice of the syzygy $\rho_1=(a_1,b_1,c_1)$ with minimal degree $d_1$, we have a morphism of graded $S$-modules
\begin{equation} \label{B1}
S(-d_1)  \xrightarrow{u} D_0(C), \  u(h)= h \cdot \rho_1.
\end{equation}
For any homogeneous syzygy $\rho=(a,b,c) \in D_0(C)_m$, consider the determinant $\Delta(\rho)=\det M(\rho)$ of the $3 \times 3$ matrix $M(\rho)$ which has as first row $x,y,z$, as second row $a_1,b_1,c_1$ and as third row $a,b,c$. Then it turns out that $\Delta(\rho)$ is divisible by $f$, see \cite{Dmax}, and we define thus a new morphism of graded $S$-modules
\begin{equation} \label{B2}
 D_0(C)  \xrightarrow{v}  S(d_1-d+1)   , \  v(r)= \Delta(\rho)/f,
\end{equation}
and a homogeneous ideal $B(C,\rho_1) \subset S$ such that $\im v=B(C,\rho_1)(d_1-d+1)$.
It is known that the ideal $B(C,\rho_1)$, when $C$ is not a free curve, defines a $0$-dimensional subscheme in $\PP^2$, which is locally a complete intersection, see \cite[Theorem 5.1]{DStJump}.

Secondly, for the reader's convenience, we recall below the main claim of \cite[Theorem 3.5]{minTjurina}, in a better formulation. Recall that 
$$t(C)=d_1+d_2-d+1$$ is the type of the curve $C$, see \cite[Definition 1.2]{NH}.

\begin{thm}
\label{thm1}
Let $C:f=0$ be a reduced curve of degree $d\geq 3$ with exponents $d_1 \leq \cdots \leq d_m$ with $m \geq 3$, and let $\rho_1$ be a non-zero syzygy of minimal degree $d_1$. Let $d'$ be the smallest integer such that $$d_3 \leq d' \leq \min(d_m,d-1)$$
and the linear system $B(C,\rho_1)_{d_1+d'-d+1}$ has a 0-dimensional base locus.
Then 
\begin{equation}
\label{eqK10}
\tau(C) \geq (d-1)^2 -d_1d_2+(d-1-d')t(C).
\end{equation}
and equality holds if and only if $C$ is a 3-syzygy curve and then $d'=d_3$.

\end{thm}

\proof
As explained in the proof of \cite[Theorem 3.5]{minTjurina}, the integer $d'$ above does indeed exist and one has
$$
\tau(C) \geq (d-1)(d-d_1-1)+d_1^2-[d_1-(d-1-d_2)][d_1-(d-1-d')].
$$
A direct simple transformation of the right hand member of the inequality in Theorem \ref{thm1} yields the new inequality
\begin{equation}
\label{eqK1}
\tau(C) \geq (d-1)^2 -d_1d_2+(d-1-d')(d_1+d_2-d+1).
\end{equation}
This yields our new formulation of  \cite[Theorem 3.5]{minTjurina}.
\endproof

By definition of $d'$ we have 
$$d-1-d' \geq 0$$
and moreover
$$t(C)=d_1+d_2-d+1 > 0$$
when $C$ is not free, by \eqref{eqF3}.
This implies the following result.

\begin{thm}
\label{thm2}
Let $C$ be a reduced plane curve of degree $d$ in $\PP^2$ with exponents $d_1 \leq d_2 \leq \ldots \leq d_m$, where $m \geq 3$.
Then $\tau(C)=(d-1)^2-d_1d_2$ if and only if  $m=3$ and $d_3=d-1$.
In all the other cases one has
$$\tau(C)>(d-1)^2-d_1d_2.$$

\end{thm}
\proof
If $m>3$, then Theorem \ref{thm1} tells us that the inequality \eqref{eqK1} is strict, and hence $\tau(C)>(d-1)^2-d_1d_2.$
When $m=3$, it follows by definition that $d'=d_3$. Then
Theorem \ref{thm1} tells us that the inequality \eqref{eqK1} is an equality
and hence we have
$$\tau(C)= (d-1)^2 -d_1d_2+(d-1-d_3)(d_1+d_2-d+1),$$
which is equivalent to the formula given in \cite[Proposition 2.1 (4)]{minTjurina} for $\tau(C)$. Therefore we get
that $\tau(C)=(d-1)^2-d_1d_2$ in this case if and only if $d_3=d-1$.
\endproof
\begin{cor}
\label{cor20}
Let $\A$ be an arrangement of $d$ lines in the projective plane $\PP^2$ with exponents $d_1 \leq d_2 \leq \ldots \leq d_m$, where $m \geq 2$. Then the following inequality holds
$$\tau(\A) \geq (d-1)^2-d_1d_2 +t(\A),$$
where $t(\A)$ is the type of the arrangement $\A$. In particular,
$\A$ is free if and only if 
$$\tau(\A)=(d-1)^2-d_1d_2 .$$ 
\end{cor}
\proof
Recall that for line arrangements the maximal degree $d_m$ in the exponents of $C$ is bounded by $d-2$, see \cite[Corollary 3.5]{Sch}. 
Then the first claim follows from the inequality \eqref{eqK10}. The second claim follows by \eqref{eqF2} and Theorem \ref{thm2}.

\endproof

\begin{cor}
\label{cor2}
Let $C$ be a reduced plane curve of degree $d$ in $\PP^2$ with exponents $d_1 \leq d_2 \leq \ldots \leq d_m$, where $m \geq 2$.
Then 
$$(d-1)^2-d_1d_2 \leq \tau(C) \leq (d-1)^2-d_1(d-1-d_1).$$ 
\end{cor}
\proof
Theorem  \ref{thm2} gives the first inequality, while the second inequality comes from the upper-bound on $\tau(C)$ given in
\cite{dPCTC}. 
\endproof

Note that $C$ is free if and only if $d_2=d-1-d_1$ and then both inequalities in Corollary \ref{cor2} become equalities.
\begin{rk}
\label{rk2}
Since $d_2 \leq d-1$, see for instance \cite[Theorem 2.4]{minTjurina}, it follows that
$$\tau(C) \geq (d-1)^2-d_1d_2 \geq (d-1)^2-d_1(d-1)=(d-1)(d-1-d_1).$$
Hence our Corollary \ref{cor2} can be regarded as an improvement of the lower bound for $\tau(C)$ given in \cite{dPCTC}.
Moreover, the equality $\tau(C)=(d-1)(d-1-d_1)$ holds if and only if
$C$ is a 3-syzygy curve with $d_2=d_3=d-1$, see \cite[Theorem 3.5]{minTjurina}
\end{rk}

\begin{ex}
\label{ex2}
Examples of 3-syzygy curves such that $d_2=d_3=d-1$ are given by the
Thom-Sebastiani curves described in \cite[Example 4.5]{minTjurina} and by the union of a smooth degree $(d-1)$ curve with one of its generic secants, see \cite[Example 4.3 (i)]{minTjurina}.

Examples of 3-syzygy curves such that $d_2=d-2$ and $d_3=d-1$ are given by the union of a smooth Fermat degree $(d-1)$ curve with a inflectional tangent meeting the curve in a single point, see \cite[Example 4.3 (ii)]{minTjurina}.

An example consisting of the Klein quartic and 4 bitangents that is a plus-one generated curve with $d_3=d-1$ can be found in 
\cite[Proposition 4.11]{JPZ}.
\end{ex}

\section{Euler polynomials of complements of plane curves} 

For any topological space $M$, having the homotopy type of a finite $CW$-complex of dimension 2, we define the Betti polynomial $B(M)$ of $M$ by the formula
\begin{equation}
\label{eqE1}
B(M)(t)=b_0(M) +b_1(M)t +b_2(M)t^2,
\end{equation}
where $b_j(M)$ denotes the $j$-th Betti number of $M$.
In particular, we can define the polynomial $B(M(C))$, where $M(C)$ is the complement $\PP^2 \setminus C$ of a reduced plane curve $C$ in $\PP^2$.
In this case,  we have the following  result.

\begin{thm}
\label{thm3}
Let $C$ be a reduced plane curve of degree $d$ in $\PP^2$ with exponents $d_1 \leq d_2 \leq \ldots \leq d_m$, where $m \geq 2$.
Let 
$$\al(C)= \tau(C)-((d-1)^2-d_1d_2).$$
If $e$ denotes the number of irreducible components of $C$, then one has 
$$(1+d_1t)(1+d_2t)-B(M(C))(t)=a(C)t+b(C)t^2,$$
where $a(C)=d_1+d_2-e+1$ and $b(C)=\mu(C)-\tau(C)+d-e+\al(C)$, with
$\mu(C)$  denoting the total Milnor  number of $C$.
In particular, one has the following.
 \begin{enumerate}
 
 \item 
 $a(C)=t(C)+(d-e) \geq 0$ and $b(C) \geq 0$;
 
 \item any of the two equalities  $a(C) = 0$ and $b(C) =0$
holds if and only if $C$ is a free line arrangement, and then the other equality also holds.

 \end{enumerate}

\end{thm}

\proof
Recall that $b_1(M(C))=e-1$ for any curve $C$.
On the other hand, the Euler numbers of the spaces under consideration satisfy
$$E(M(C))=E(\PP^2) -E(C)=3-E(C) \text{ and } E(C)=E(C_d)+\mu(C),$$
where $C_d$ is a smooth curve of degree $d$. Hence
$$E(C_d)=2-2g=2-(d-1)(d-2).$$
Putting these formula together yields 
$$b_2(M(C))=(d-1)^2-\mu(C)-(d-e),$$
and this clearly prove our first claim. For the second claim we notice that
in general we have
$$d_1+d_2 \geq d-1 \geq e-1$$
which implies that $a(C) \geq 0$. The equality $a(C)=0$ implies
 $t(C)=d_1+d_2-d+1=0$ and $e=d$, namely that $C$ is a free line arrangement. 
 
 Moreover, $\mu(C) \geq \tau(C)$, with equality if and only if all the singularities of $C$ are quasi homogeneous. Corollary \ref{cor2} shows that $\al(C) \geq 0$ and since $d \geq e$ obviously, we get $b(C) \geq 0$. Finally, $b(C)=0$ implies that $d=e$, and hence $C$ is a line arrangement, and also yields $\al(C)=0$.  Corollary \ref{cor20} implies that $C$ is free 
 and this completes the proof of the second claim.
\endproof

\begin{cor}
\label{cor3}
Let $C$ be a reduced plane curve of degree $d$, having $e$ irreducible components. 
If $C$ is free with exponents $(d_1,d_2)$, then
$$(1+d_1t)(1+d_2t)=B(M(C))(t)+(d-e)t(1+t)+(\mu(C)-\tau(C))t^2.$$
In particular, if in addition all the singularities of $C$ are quasi homogeneous, one has
$$(1+d_1t)(1+d_2t)=B(M(C))(t)+(d-e)t(1+t)$$
and the Euler number $E(M(C))$ satisfies
$$E(M(C))=B(M(C))(-1)=(d_1-1)(d_2-1).$$
\end{cor}
\proof
For a free curve $C$ we have $d_1+d_2=d-1$ and $\al(C) =0$,
recall \eqref{eqF1} and \eqref{eqF2}.
\endproof
We notice that the equality \eqref{eqF4} follows also from Corollary \ref{cor3}, since $e=d$ in this case,
see also \cite{Pok}.

\begin{rk}
\label{rk1}
The reduced plane curves with $d_1=1$ were studied in \cite{dPCTC2}, where the interested reader may find out many explicit examples.
It turns out that the  curves with $d_1=1$ are either free or nearly free, that is 3-syzygy curves such that  the exponents are $(1,d-1,d-1)$. This follows from  \cite[Proposition 1.3 (ii)]{dPCTC2} and the characterization of nearly free curves given in \cite{Dmax}. Moreover, such a curve admits a 1-dimensional connected group of symmetries $H$, see \cite{dPCTC2}. Since
$H$ is either $(\C^*, \cdot)$ or $(\C,+)$, it follows that the closures of the 1-dimensional orbits of $H$ are rational curves. This implies that any irreducible component of a plane curve $C$ with $d_1=1$ is rational.
When $C$ is free, we can give the following alternative proof for this result.
\end{rk}

\begin{cor}
\label{cor31}
Let $C$ be a free curve with exponents $(d_1,d_2)$ such that $d_1=1$. Then any irreducible component of $C$ is a rational curve. 
\end{cor}
\proof
Using Corollary \ref{cor3} we get 
$$E(C)=B(M(C))(-1)=(d_1-1)(d_2-1)-(\mu(C)-\tau(C)) \leq 0,$$
since $\mu(C) \geq \tau(C)$.
Then using a conjecture of W. Veys in \cite{W}, proved by  A. de Jong and  J. Steenbrink in \cite{JS} and by R. Gurjar and  A. Parameswaran in \cite{GP}, it follows that
any irreducible component of $C$ is rational.
\endproof

\section{Two results about line arrangements}

First we recall two main results from \cite{POG} in the setting of line arrangements in $\PP^2$. Similar results appear also in
\cite[Theorem 1.11]{A} and in \cite[Theorems 3.5 and 3.6]{MP}. Since all the singularities of line arrangements are quasi-homogeneous, it follows that all the invariants $\epsilon$ in the quoted results are $0$.
First we restate \cite[Theorem 1.3]{POG}.
\begin{thm}
\label{thm1B}
Let $\A'$ be a  line arrangement in $\PP^2$, $L$ a line in $\PP^2$, which is not in $\A'$. We assume that the union $\A=\A' \cup L$ is a free curve with exponents $(d_1,d_2)$. Then  the exponents $(d_1',d_2')$ (resp. $(d_1',d_2',d_3')$) of the free (resp. plus-one generated) arrangement $\A'$ and  the number $r=|\A' \cap L|$ of intersection points satisfy one of the following conditions, and all these three cases are possible.
\begin{enumerate}

\item $d_1<d_2$, $d_1'=d_1$ and $d_2'=d_2-1$. In this case $\A'$ is  free  and $r=d_1+1$.

\item  $d_1'=d_1-1$ and $d_2'=d_2$. In this case $\A'$ is  free  and $r=d_2+1$.

\item $d_1'=d_1$ and $d_2'=d_2$. In this case $\A'$ is  plus-one generated  and $$r=|\A'|-d_3' \leq |\A|-1-d_2=d_1.$$

\end{enumerate}

In particular, $\A'$ is  free  if and only if $r\geq d_1+1$.
\end{thm}

Now we restate \cite[Theorem 1.4]{POG}.

\begin{thm}
\label{thm2B}
Let $\A'$ be a  line arrangement in $\PP^2$, $L$ a line in $\PP^2$, which is not in $\A'$.
We consider the union $\A=\A' \cup L$ and assume that $\A'$ is a free arrangement with exponents $(d_1',d_2')$. Then  the exponents $(d_1,d_2)$ (resp. $(d_1,d_2,d_3)$) of the free (resp. plus-one generated)  $\A$ and the number $r=|\A' \cap L|$ of intersection points  satisfy one of the following conditions, and all these three cases are possible.
\begin{enumerate}

\item $d_1=d_1'$ and $d_2=d_2'+1$. In this case $\A$ is  free  and $r=d_1'+1$.

\item $d_1'<d_2'$, $d_1=d_1'+1$ and $d_2=d_2'$. In this case $\A$ is free  and $r=d_2'+1$.

\item $d_1=d_1'+1$ and $d_2=d_2'+1$. In this case $\A$ is  plus-one generated  and $$r=d_3+1 \geq d_2+1 =d_2'+2.$$

\end{enumerate}
In particular, $\A$ is  free  if and only if $r\leq d_2'+1$.

\end{thm}

We recall that (a partial case of)  Terao's conjecture says that if
$\A$ and $\B$ are two line arrangements in $\PP^2$ having the same combinatorics, and if $\A$ is free, then $\B$ is also free.
By a {\it minimal counter example to Terao's Conjecture} we mean a pair $\A, \B$ as above such that Terao's conjecture holds for all pairs
$\A',\B'$ with $|\A'|=|\B'| <|\A|=|\B|$, but not for the pair $\A, \B$.
Our next result  gives some information on such a  minimal counter example to Terao's Conjecture, supposing it exists. A slightly stronger version of this result follows from \cite{Ab}, see Remark \ref{rk40} below for details.

\begin{thm}
\label{thm4}
Let $\A, \B$ be a  minimal counter example to Terao's Conjecture, such that $\A$ is free with exponents $d_1 \leq d_2$.
Then the following  property holds for the free arrangement $\A$:
  for any line $L \in \A$, the number $r_L$ of multiple points of $\A$
 situated on $L$ satisfies the inequality 
 $$r_L  \leq d_1.$$

\end{thm}
\proof
We prove that if there is a line $L \in \A$ such that
$$r_L \geq d_1+1$$
 then $\B$ is also free, and hence $\A,\B$ is not a counter example to Terao's Conjecture.
If such a line $L$ exists, then we denote by $L'$ the line in $\B$ corresponding to the line $L \in \A$ with $r_L  \geq d_1+1$ under the isomorphism of intersection lattices $L(\A) \simeq L(\B)$. 
Let
$\A'=\A \setminus \{L\}$ and $\B'=\B \setminus \{L'\}$ be the two deleted arrangements obtained from $\A$ and $\B$ by deleting the line $L$ and respectively $L'$. Clearly
$$r_{L‘}=|\B' \cap L'| =r_L.$$
We apply Theorem \ref{thm1B} and conclude that
$\A'$ is a free line arrangement. Indeed, this happens in the cases $(1)$ and $(2)$ of  Theorem \ref{thm1B}. In the remaining case $(3)$,  we have
$r_L=|\A' \cap L| \leq d_1,$
hence this case cannot occur in our situation.

On the other hand, it is clear that the intersection lattice isomorphism $L(\A) \simeq L(\B)$ yields a new intersection lattice isomorphism $L(\A') \simeq L(\B')$.
Since $\A,\B$ was supposed to be a  {\it minimal } counter example to Terao's Conjecture, it follows that Terao's Conjecture holds for the pair
$\A',\B'$ and hence $\B'$ is also a free line arrangement with
the same exponents as $\A'$, call them $d_1' \leq d_2'$.
We apply now Theorem \ref{thm2B} to the arrangements $\B'$ and $\B$, and note that in case (1) of 
Theorem \ref{thm1B} we have
$$r_{L'}=r_L=d_1'+1=d_1+1 \leq d_2+1$$
and hence $\B$ is free by Theorem \ref{thm2B}.
Similarly,  in case (2) of 
Theorem \ref{thm1B} we have
$$r_{L'}=r_L=d_1+1=d_2+1=d_2'+1$$
and hence again $\B$ is free.

\endproof

\begin{rk}
\label{rk40}
The weaker inequality $r_L  \leq d_1+1$ holds in fact for any free line arrangement
$\A$ with exponents $d_1 \leq d_2$, unless one has $r_L=d_2+1$, see \cite[Corollary 1.2]{Ab}. There is a difference in our notations, since in 
\cite{Ab} affine line arrangements in $\C^2$ are considered.
On the other hand, it is clear by \cite[Theorem 1.1 (3)]{Ab} that a free arrangement $\A$ with $r_L=d_1+1$ or $r_L=d_2+1$ cannot be a counter example to Terao's Conjecture. Hence in fact any counter example to Terao's Conjecture has to satisfy $r_L \leq d_1$.
In relation to these facts, see also
\cite[Proposition 5.2 and Proposition 5.3]{FV}.
\end{rk}

\begin{rk}
\label{rk4}
The monomial line arrangement
$$\A=\A(m,m,3): (x^m-y^m)(y^m-z^m)(x^m-z^m)=0$$
is free with exponents $d_1=m+1$ and $d_2=2m-2$ for $m \geq 3$, see for instance \cite[Example 8.6 (i)]{DHA}. For any line $L \in \A$ one has
$r_L=m+1$ so the condition $r_L \leq d_1$ from Theorem \ref{thm4} holds. On the other hand,  the full monomial line arrangement
$$\A=\A(m,1,3): xyz(x^m-y^m)(y^m-z^m)(x^m-z^m)=0$$
is free with exponents $d_1=m+1$ and $d_2=2m+1$ for $m \geq 2$, see for instance \cite[Example 8.6 (ii)]{DHA}. For any line $L \in \A$ one has
$r_L=m+2$ so the condition $r_L \leq d_1$ from Theorem \ref{thm4} does not hold. Similarly, this condition fails for the Hessian line arrangement which is free with exponents $d_1=4$ and $d_2=7$ and
for any line $L $ in it one has
$r_L=5$,   see for instance \cite[Example 8.6 (i)]{DHA}. We believe that the condition
$r_L \leq  d_1$ from Theorem \ref{thm4} limits a lot the possibilities for $\A$, but we have no result
in this direction for the moment.

\end{rk}

Let $\A$ be a line arrangement with exponents $d_1 \leq \ldots \leq d_m$ and let $m(\A)$ be the maximal multiplicity of a multiple point in $\A$.
Then Lemma 5.2 in \cite{Sch} states (without a proof and in a rather cryptic way) that for the line arrangement $\A$ one has
\begin{equation}
\label{eqSS}
 m(\A)-1 \leq d_m
\end{equation}
see Theorem 5.4 claim 5. in \cite{Sch} for a clearer statement, where $M=m(\A)-1$. In fact this claim can be improved as follows. 

\begin{thm}
\label{thm5}
Let $m(\A)$ be the maximal multiplicity of a point in the line arrangement $\A$,
and $n(\A)$ the maximal multiplicity of a point in $\A \setminus \{p\}$, where $p$ is any point in $\A$ of multiplicity $m(\A)$.
Then one has  
$$ m(\A)-1 \leq    t(\A) +m(\A)-1 \leq  d_2  \leq d-n(\A),$$
where $t(\A)=d_1+d_2-d+1 \geq 0$ is the type of the arrangement $\A$ and $d$ is the number of lines in $\A$.
In particular, the equality $d_2=m(\A)-1$ implies that $\A$ is free and $d_1=d-m(\A)$, and for a non free line arrangement $\A$ one has 
$$ m(\A)\leq d_2.$$
Moreover, the type $t(\A)$ of the line arrangement $\A$ satisfies
$$0 \leq t(\A)\leq d+1-m(\A)-n(\A).$$
\end{thm}
\proof

The inequality
$$ (d_1+d_2-d+1) + m(\A)-1 \leq d_2$$
is clearly equivalent to the inequality
\begin{equation}
\label{eqSS1}
d_1 \leq d-m(\A),
\end{equation}
which is proven in \cite[Theorem 1.2]{Mich}. Hence the lower bound for $d_2$ is obtained. To get the upper bound, we use \cite[Lemma 3.2]{minTline}, which implies that we have 
\begin{equation}
\label{eqSS2}
d_2 \leq d-n(\A).
\end{equation}
More precisely, let $p,q \in \A$ be two distinct points such that the multiplicity of $\A$ at $p$ (resp. $q$) is $m=m(\A)$ (resp. $n=n(\A)$). Then
there are two primitive syzygies $\rho_p \in D_0(\A)_{d-m}$ and
$\rho_q \in D_0(\A)_{d-n}$, see \cite[Lemma 3.2]{minTline}. It follows that
$$d_1 \leq d-m \leq d-n.$$
If $d_2 >d-n$ and $\rho$ is a generator for $D_0(\A)_{d_1}$, then there are two nonzero polynomials $P,Q \in S$ such that
$$\rho_p=P\rho \text{ and } \rho_q=Q\rho.$$
It follows that 
\begin{equation}
\label{eqSS3}
\rho_q=A\rho_p
\end{equation}
where $A=Q/P$ is a rational function. It was shown in  \cite[Lemma 3.2]{minTline} that the equality \eqref{eqSS3} is impossible with $A$ a polynomial. Exactly the same proof shows that the equality \eqref{eqSS3} is impossible with $A$ a rational function. This implies that
$d_2 \leq d-n$ as we have claimed.
These two inequalities  \eqref{eqSS1} and \eqref{eqSS2} prove all the remaining claims.
\endproof

\begin{rk}
\label{rk5}
Consider the line arrangement
$$\A: f=(x^5-y^5)(x+2y+z)(x+3y-5z)=0$$
consisting of 5 lines through a point plus two generic lines. A direct computation shows that $\A$ has exponents $(2,5,5)$, hence $\A$ is a nearly free arrangement satisfying 
$$5=m(\A)=d_2=d-n(\A)= 7-2.$$
It follows that the inequalities in Theorem \ref{thm5} are sharp.
The corresponding Poincaré polynomial is
$$\pi(\A,t)=(1+3t)^2.$$
This is a special case of Example 5.3 in \cite{Sch}, which seem to suggest that in that paper Lemma 5.2 {\it refers only to free arrangements}. Indeed, the author says that we can use this result in order to show that $\A$ is not free, perhaps thinking that the Poincaré polynomial yields the exponents as in the free case.
\end{rk}

\begin{rk}
\label{rk6}
There is the following relation between Theorems \ref{thm4} and \ref{thm5} above. Consider the arrangement $\A$ in Theorem \ref{thm4}
and let $p \in \A$ be a point of maximal multiplicity $m(\A)$. If $L$ is a line in $\A$ not passing through $p$, it is clear that
$$r_L \geq m(\A)$$
since any line through $p$ intersects $L$ in a distinct point.
If this inequality is an equality for any line $L$ not passing through $p$,
then $p$ is a modular point, and hence $\A$ is a supersolvable line arrangement. Since for  supersolvable line arrangements Terao's conjecture holds, it follows that this is not the case, and hence there is a line $L$ not through $p$ with $r_L >m(\A)$. This yields the stronger inequality 
$$d_1 > m(\A),$$
for the arrangement $\A$ in Theorem \ref{thm4}, which also follows  from
\cite[Corollary 1.4]{Mich}. This result  says that for free arrangements $\A$ with $m(\A) \geq d_1$ Terao's Conjecture holds.
 
\end{rk}


\end{document}